\documentclass{amsproc}
\usepackage{mathrsfs}
\usepackage[OT2,T1]{fontenc}
\usepackage{hyperref}
\usepackage{bm}

\newcommand{\comment}[1]{}
\newcommand{\ms}[1]{\mathscr{#1}}
\newcommand{\mc}[1]{\mathcal{#1}}

\newcommand{\RR}{{\mathbb R}}
\newcommand{\CC}{{\mathbb C}}
\newcommand{\PP}{{\mathbb P}}
\newcommand{\ZZ}{{\mathbb Z}}
\newcommand{\QQ}{{\mathbb Q}}
\newcommand{\HH}{{\mathbb H}}
\newcommand{\FF}{{\mathbb F}}

\newcommand{\MM}{\mathbb{M}}
\newcommand{\Th}{{\operatorname{Th}}}    
\newcommand{\ON}{{\operatorname{ON}}}    
\newcommand{\SL}{\operatorname{\textsl{SL}}}      
\newcommand{\GL}{\operatorname{\textsl{GL}}}      
\newcommand{\Aut}{\operatorname{Aut}}
\newcommand{\Sel}{\operatorname{Sel}}

\newcommand{\TP}{\operatorname{\textsl{P}}\!}
\newcommand{\TS}{\operatorname{\textsl{T}}\!}
\DeclareSymbolFont{cyrletters}{OT2}{wncyr}{m}{n}\DeclareMathSymbol{\Sha}{\mathalpha}{cyrletters}{"58}
\newcommand{\tr}{\operatorname{{tr}}}
\newcommand{\xmod}{{\rm \;mod\;}}
\newcommand{\tor}{{\rm tor}}

\copyrightinfo{2009}{American Mathematical Society}

\newtheorem{theorem}{Theorem}[section]

\theoremstyle{definition}

\theoremstyle{remark}

\numberwithin{equation}{section}

\begin{document}

\title{From the Monster to Thompson to O'Nan}

\author{John F. R. Duncan}

\address{Department of Mathematics, Emory University, Atlanta, GA 30322, U.S.A.}

\email{john.duncan@emory.edu}

\dedicatory{To Geoff Mason, from whom we have learnt much about moonshine.}

\thanks{We thank Jia-Chen Fu, Maryam Khaqan, Michael Mertens, Ken Ono and an anonymous referee for helpful comments on earlier versions of this work, and we gratefully acknowledge financial support from the U.S. National Science Foundation (DMS 1601306).}


\date{}

\begin{abstract}
The commencement of monstrous moonshine is a connection between the largest sporadic simple group---the monster---and complex elliptic curves. Here we explain how a closer look at this connection leads, via the Thompson group, to recently observed 
relationships between the non-monstrous sporadic simple group of O'Nan and 
certain families of elliptic curves defined over the rationals. We also describe umbral moonshine from this perspective.
\end{abstract}

\maketitle

\section{Elliptic Curves}\label{sec:ellipticcurves}

Since they are the driving force in this work we begin by recalling that
an {\em elliptic curve} over a field $K$ is a pair $(E,O)$ where $E$ is a non-singular projective algebraic curve over $K$ with genus one and $O$ is a point of $E$ defined over $K$. 
More concretely, and assuming for simplicity that the characteristic of $K$ is not $2$ or $3$, any elliptic curve over $K$ is isomorphic to one of the form $(E,O)$ where $E\subset \PP^2(K)$ is specified by a homogeneous cubic 
\begin{gather}\label{eqn:ellipticcurvehomogeneousequation}
Y^2Z=X^3+AXZ^2+BZ^3
\end{gather} 
for $A,B\in K$ such that the {\em discriminant} $\Delta:=-16(4A^3+27B^2)$ is not zero (cf. e.g. \cite{MR1193029,MR2514094}). The distinguished point is $O=[0,1,0]$. 
There is no loss of information in suppressing the powers of $Z$ in (\ref{eqn:ellipticcurvehomogeneousequation}) so from now on we do this and write simply 
\begin{gather}\label{eqn:ellipticcurveequation}
E\,:\,
y^2=x^3+Ax+B
\end{gather} 
for the curve in $\PP^2(K)$ defined by (\ref{eqn:ellipticcurvehomogeneousequation}).

By force of the genus condition 
the set $E(K)$ of {\em $K$-rational points} (i.e. points defined over $K$) naturally acquires an abelian group structure with identity element $O$.
A fundamental result about $E(K)$ is that it is finitely generated as an abelian group if $K$ is a {\em number field} (i.e. a finite extension of the rationals, $\QQ$). 
This was proved for $K=\QQ$ 
by Mordell \cite{Mor_RtnSlnIndEqn3rd4thDgr}, and subsequently generalized to number fields (and higher dimensional abelian varieties) by Weil \cite{MR3532958}. 

The computation of the rank of $E(K)$, for $K$ a number field, is a challenging problem (cf. e.g. \cite{MR2238272}). 
To get a taste for it let $E_{15}'$ and $E_{15}''$ be the elliptic curves 
over $\QQ$ defined by 
\begin{gather}
\begin{split}\label{eqn:E1E2}
	E_{15}'\,:\, 
	y^2&=x^3-831168x + 134894592,\\ 
	E_{15}''\,:\, 
	y^2&=x^3-60051888x + 82842141312.
\end{split}
\end{gather} 
Then, as we will see in \S\ref{sec:ECONan} (wherein the significance of the subscripts will also be revealed), the group $E_{15}'(\QQ)$ is finite but $E_{15}''(\QQ)$ is infinite. Perhaps neither of these facts is particularly obvious from the defining expressions (\ref{eqn:E1E2}).

The {\em $j$-invariant} of $E$ as in (\ref{eqn:ellipticcurveequation}) is 
\begin{gather}
j(E):=1728\frac{4A^3}{4A^3+27B^2}.
\end{gather}
With this definition it develops (see e.g. \S1 in Chapter III of \cite{MR2514094}) that elliptic curves $E$ and $E'$ over $K$ are isomorphic over an algebraic closure $\overline K$ of $K$ if and only if $j(E)=j(E')$. Also, for any $j_0\in \overline K$ there exists an elliptic curve $E_0$, defined over $K(j_0)$, such that $j(E_0)=j_0$.
So in particular, taking $K=\CC$ to be the complex numbers we find that $j$ defines a bijection between isomorphism classes of complex elliptic curves and $\CC$.

For complex elliptic curves the {\em upper half-plane} $\HH:=\{\tau\in \CC\mid \Im(\tau)>0\}$ plays a special role. 
This is because for any $\tau\in \HH$ the quotient 
\begin{gather}
E_\tau:=\CC/(\ZZ+\ZZ\tau)
\end{gather}
defines an elliptic curve over $\CC$ (with $O$ represented by $0$), and every complex elliptic curve is isomorphic to $E_\tau$ for some $\tau\in \HH$. 
The $j$-invariant of $E_\tau$ satisfies 
\begin{gather}\label{eqn:jEtau}
	j(E_\tau)=
		\frac{\eta(\tau)^{24}}{\eta(2\tau)^{24}}+4096\frac{\eta(\tau)^{24}}{\eta(\frac12\tau)^{24}}-4096\frac{\eta(\frac12\tau)^{24}\eta(2\tau)^{24}}{\eta(\tau)^{48}}
		+768
\end{gather}
where 
$\eta(\tau):=e^{\pi i\tau\frac1{12}}\prod_{n>0}(1-e^{2\pi i\tau n})$ is the {\em Dedekind eta function}.
So the $j$-invariant defines a holomorphic function on $\HH$. 
We have $E_\tau\simeq E_{\tau'}$ if and only if 
$\gamma\tau=\tau'$ for some element $\gamma$ in the {\em modular group} $\SL_2(\ZZ)$, where the action 
is 
\begin{gather}\label{eqn:SL2ZactiononHH}
\left(\begin{matrix} a&b\\c&d\end{matrix}\right)\tau:=\frac{a\tau+b}{c\tau+d}.
\end{gather} 
So $j(\tau):=j(E_\tau)$ defines correspondences between the set of isomorphism classes of complex elliptic curves, the orbit space $\SL_2(\ZZ)\backslash\HH$, and 
$\CC$.

Set 
$q:=e^{2\pi i\tau}$. Using (\ref{eqn:jEtau}) and the formula $\eta(\tau)=q^{\frac1{24}}\prod_{n>0}(1-q^n)$ we may directly compute the first few terms in the {\em Fourier series} expansion
\begin{gather}
\label{eqn:jtauqseries}
	j(\tau)=
		q^{-1}+744+196884q+21493760q^2+864299970q^3+\dots
\end{gather}

\section{Monstrous Moonshine}\label{sec:monstrousmoonshine}

We now turn to the {\em Fischer--Griess monster} $\MM$ which is the finite simple group with order 
\begin{gather}\label{eqn:orderMM}
	\#\MM= 808017424794512875886459904961710757005754368000000000
\end{gather}
that was discovered independently by Fischer and Griess (cf. \cite{MR0399248,MR671653}) 
in $1973$, and is distinguished as the largest of the sporadic simple groups (cf. \cite{atlas}).

It was some time after its discovery that the existence of the monster group was verified by Griess \cite{MR671653}. In the interim Conway--Norton conjectured that there is an embedding of $\MM$ into $\GL_n(\CC)$ for $n=196883$, and McKay observed (cf. \cite{MR554399,Tho_NmrlgyMonsEllModFn}) that 
this number is just $1$ away from being the coefficient of $q$ in the Fourier series expansion (\ref{eqn:jtauqseries}) of the $j$-invariant for complex elliptic curves.

Thompson extended McKay's observation \cite{Tho_NmrlgyMonsEllModFn} using the (conjectural) character table for $\MM$ that had been computed by this point by Fischer--Livingstone--Thorne (cf. \cite{MR554399,atlas}). Inspired by these coincidences he suggested the existence of a faithful graded infinite-dimensional $\MM$-module 
\begin{gather}\label{eqn:MMmodule}
V=\bigoplus_n V_n
\end{gather} 
with $j(\tau)-744=\sum_n \dim V_n q^n$. He also suggested to consider the graded trace functions 
\begin{gather}\label{eqn:MTseries}
T_g(\tau):=\sum_n \tr(g|V_n)q^n
\end{gather} 
for $g\in \MM$, now known as {\em McKay--Thompson series} for $\MM$, and proposed \cite{Tho_FinGpsModFns} that each such function is a {principal modulus} for a genus zero subgroup of $\SL_2(\RR)$. 

We recall now that if $\Gamma<\SL_2(\RR)$ is {\em commensurable} with the {modular group} $\SL_2(\ZZ)$, in the sense that $\Gamma\cap\SL_2(\ZZ)$ has finite index in both $\Gamma$ and $\SL_2(\ZZ)$, then the action of $\Gamma$ on $\HH$ extends naturally to $\PP^1(\QQ)=\QQ\cup\{\infty\}$ (cf. e.g. Proposition 2.13 in \cite{Dun_ArthGrpsAffE8Dyn}), and the
quotient 
\begin{gather}\label{eqn:XGamma}
X_\Gamma:=\Gamma\backslash(\HH\cup\PP^1(\QQ))
\end{gather}
is naturally a compact complex curve (cf. e.g. \cite{Shi_IntThyAutFns}). Such a group $\Gamma$ is called {\em genus zero} if $X_\Gamma\simeq \PP^1(\CC)$, i.e., if $X_\Gamma$  has genus zero as a real surface.
A $\Gamma$-invariant holomorphic function on $\HH$ is called a {\em principal modulus} or {\em Hauptmodul} for $\Gamma$ if it extends from $\Gamma\backslash\HH\to \CC$ to an isomorphism $X_\Gamma\xrightarrow{\sim}\PP^1(\CC)$. 

We will only consider subgroups $\Gamma<\SL_2(\RR)$ that are commensurable (in the above sense) with $\SL_2(\ZZ)$ in this work. The {\em cusps} of such a group are the points of 
\begin{gather}\label{eqn:cusps}
\Gamma\backslash\PP^1(\QQ)\subset X_\Gamma. 
\end{gather}
The {\em infinite cusp} is the cusp represented by $\infty\in \PP^1(\QQ)$, and we call any other cusp {\em non-infinite}.

Conway--Norton computed explicit candidates for the $T_g$ and made numerous observations regarding them in \cite{MR554399}. 
From this emerged the {\em monstrous moonshine conjecture}, being the statement that 
the monster module (\ref{eqn:MMmodule})
is such that 
for each $g\in \MM$ the corresponding McKay--Thompson series 
(\ref{eqn:MTseries}) 
is the {normalized principal modulus} for a {genus zero group} 
\begin{gather}\label{eqn:Gammag}
\Gamma_g<\SL_2(\RR),
\end{gather} 
commensurable with $\SL_2(\ZZ)$,
specified explicitly in \cite{MR554399}. 

A principal modulus $T$ for a genus zero group $\Gamma<\SL_2(\RR)$ is said to be {\em normalized} if 
\begin{gather}\label{eqn:normalizedpm}
T(\tau)=
\begin{cases}
q^{-1}+O(q)& \text{ as $\Im(\tau)\to \infty$,} \\
O(1)&\text{ as $\tau$ tends to any non-infinite cusp of $\Gamma$}.
\end{cases}
\end{gather}
The significance of this is that a normalized principal modulus is uniquely determined by its invariance group $\Gamma$, because if $T$ and $T'$ are two normalized principal moduli for $\Gamma$ then the difference $T-T'$ defines a bounded holomorphic function on $X_\Gamma$ (\ref{eqn:XGamma}). The only such functions are constants, and $T(\tau)-T'(\tau)=O(q)$ 
as $\Im(\tau)$ goes to $\infty$ by (\ref{eqn:normalizedpm}), so $T-T'$ must vanish identically.

Because principal moduli exist only for genus zero groups the fact that the monstrous McKay--Thompson series are normalized principal moduli is sometimes referred to as the {\em genus zero property} of monstrous moonshine. It is this property that gives monstrous moonshine its power, because it allows to compute the structure of $V$ (\ref{eqn:MMmodule}) as a module for $\MM$ (\ref{eqn:orderMM}) with no more information than the assignment $g\mapsto \Gamma_g$ (\ref{eqn:Gammag}) of Conway--Norton.

The function $j(\tau)-744$ is the normalized principal modulus for the {modular group}, 
so 
\begin{gather}
X_{\SL_2(\ZZ)}=\SL_2(\ZZ)\backslash (\HH\cup \PP^1(\QQ))
\end{gather} 
is a (coarse) moduli space for complex elliptic curves.
Conway--Norton found that for each $g\in \MM$ the group $\Gamma_g$ (\ref{eqn:Gammag}) normalizes a {\em Hecke congruence subgroup}
\begin{gather}\label{eqn:Gamma0N}
\Gamma_0(N_g):=\left\{\left(\begin{matrix}a&b\\c&d\end{matrix}\right)\in\SL_2(\ZZ)\mid c\equiv 0 \xmod N_g\right\}
\end{gather} 
for some $N_g$. For this reason the curves $X_{\Gamma_g}$ for $g\in \MM$ (cf. (\ref{eqn:XGamma}), (\ref{eqn:Gammag})) may be interpreted as (coarse) moduli spaces of {\em isogenies} (i.e. homomorphisms) of complex elliptic curves with cyclic kernel (cf. e.g. \S6 in Chapter IX of \cite{MR1193029}). 

The monstrous moonshine conjecture is a theorem thanks to work of Frenkel--Lepowsky--Meurman \cite{FLMPNAS,FLMBerk,FLM} and Borcherds \cite{Bor_PNAS,MR1172696}. See \cite{mnstmlts,MR2201600} for reviews of this. Curiously, their results do not so far seem to have told us anything new about elliptic curves or their isogenies, other than that they are connected, as above, to the monster.

\section{Thompson Moonshine}\label{sec:thompsonmoonshine}

We mentioned in \S\ref{sec:ellipticcurves} that the computation of $E(K)$ for $K$ a number field is a prominent problem in elliptic curves. 
At first glance $j(\tau)$ (\ref{eqn:jtauqseries}) does not seem to have anything to say about this because by passing to $E(\CC)$ we ``wash away'' the subtle arithmetic of $K$. For example, $E_{15}'$ and $E_{15}''$ as in (\ref{eqn:E1E2}) have the same $j$-invariant, $j(E_{15}')=j(E_{15}'')=
\frac{111284641}{50625}$, 
so they are the same over $\CC$, even though $E_{15}'(\QQ)$ is finite and $E_{15}''(\QQ)$ is infinite (cf. \S\ref{sec:ECONan}). However, hints of richer structure emerge when we consider the values of $j(\tau)$ at certain special points in $\HH$.

Say that $\tau\in \HH$ is a {\em CM point} if $\QQ(\tau)$ is a quadratic extension of $\QQ$. 
Then it develops (see e.g. \cite{MR0201394}) that $j(\tau)$ is an algebraic integer when $\tau$ is CM. 
An algebraic integer written in the form $j(\tau)$ for $\tau$ a CM point 
is called a {\em singular modulus}. 
For example, 
we have 
\begin{gather}\label{eqn:jtau85995}
j\left(\frac{1+\sqrt{-15}}2\right)=-52515-85995\frac{1+\sqrt{5}}{2}.
\end{gather}

It turns out that this identity (\ref{eqn:jtau85995}) has significance for sporadic simple groups. 
To explain this we note that the monster has a conjugacy class, called $3C$, such that $T_{g}$ satisfies $T_{g}(\tau)^3=j(3\tau)$ when $g\in 3C$. 
In this case the centralizer $C_\MM(g)$ takes the form 
\begin{gather}
C_\MM(3C)\simeq
\ZZ/3\ZZ\times \Th
\end{gather} 
(cf. \cite{atlas}) 
where $\Th$ denotes the sporadic simple {\em Thompson group} \cite{MR0409630,MR0399193}, having order
\begin{gather}\label{eqn:orderTh}
\#\Th= 90745943887872000.
\end{gather}
The 
connection to (\ref{eqn:jtau85995}) is that 
the Thompson group
has a pair of irreducible representations of dimension $85995$ whose character values lie in $\QQ(\sqrt{-15})$ (cf. \cite{atlas}).

This observation may serve as a starting point for the {\em Thompson moonshine} introduced by Harvey--Rayhaun \cite{Harvey:2015mca}, which associates a 
(weakly holomorphic) modular form
\begin{gather}\label{eqn:ThMTseries}
\check F_g^{\Th}(\tau)=2q^{-3}+\sum_{\substack{D\geq 0\\D\equiv 0,1\xmod 4}}C^\Th_g(D)q^{D}
\end{gather} 
of weight $\frac12$ for $\Gamma_0(4o(g))$ to each $g\in \Th$, in such a way that the function $g\mapsto (-1)^DC^\Th_g(D)$ is 
a character of $\Th$ for each fixed $D\geq 0$. (This latter fact was confirmed by Griffin--Mertens in \cite{MR3582425}.)
For $g=e$ the identity element we have
\begin{gather}\label{eqn:Thgdim}
\check F_e^{\Th}(\tau)=2q^{-3}+248+54000q^4-171990q^5+\dots,
\end{gather} 
and, as explained in \S2 of \cite{Harvey:2015mca},  
the coefficients $C_e^\Th(D)$ may be expressed 
as linear combinations of
singular moduli. For example, we have
$C_e^\Th(5)=-2\times 85995$, and 
\begin{gather}\label{eqn:85995asatrace}
	85995 = \frac{1}{\sqrt{5}}\left(
	j\left(\frac{1+\sqrt{-15}}4\right)
	-
	j\left(\frac{1+\sqrt{-15}}2\right)
	\right)
\end{gather}
is a weighted sum of 
singular values of $j$ 
at 
CM points 
involving $\sqrt{-15}$. (In particular $\frac{1+\sqrt{-15}}4$ is also a CM point. We will describe a general framework for such sums in \S\ref{sec:tracesofsingularmoduli}.)

The McKay--Thompson series (\ref{eqn:ThMTseries}) of Thompson moonshine 
share a property analogous to the normalized principal modulus property (\ref{eqn:normalizedpm}) that characterizes the McKay--Thompson series (\ref{eqn:MTseries}) of monstrous moonshine. To formulate this cleanly it is best to pass to the vector-valued functions $F^\Th_g:=(F^\Th_{g,0},F^\Th_{g,1})$ where
\begin{gather}\label{eqn:vvFTH}
	F^\Th_{g,r}(\tau):=\delta_{r,1}2q^{-\frac34}+\sum_{\substack{D\geq 0\\D\equiv r^2\xmod 4}} C^\Th_g(D)q^{\frac{D}4}
\end{gather}
for $r\in\{0,1\}$. Then $F^\Th_g$ is a (weakly holomorphic) vector-valued modular form of weight $\frac12$ for $\Gamma_0(o(g))$ 
(\ref{eqn:Gamma0N}), whereas $\check F^\Th_g$ (\ref{eqn:ThMTseries}) is only modular for the subgroup $\Gamma_0(4o(g))$.
Also, for each $g\in \Th$ the corresponding vector-valued McKay--Thompson series $F_g^\Th$ is uniquely determined---up to a theta series---amongst weakly holomorphic modular forms 
of weight $\frac12$ (with a suitably defined multiplier) for $\Gamma_g^{\Th}:=\Gamma_0(o(g))$ by the property that 
\begin{gather}\label{eqn:Thpmanalogue}
F_g^\Th(\tau)=
\begin{cases}
(0,2q^{-\frac34})+O(1) &\text{ as $\Im(\tau)\to \infty$,}\\ 
O(1) &\text{ as $\tau$ tends to any non-infinite cusp of $\Gamma_g^{\Th}$.}
\end{cases}
\end{gather}

So although the groups $\Gamma_g^{\Th}$ are generally not genus zero (e.g. we have $\Gamma_g^\Th=\Gamma_0(31)$ when $o(g)=31$, and the genus of $X_{\Gamma_0(31)}$ (\ref{eqn:XGamma}) is $2$) we still have similar predictive power over the $F_g^\Th$ to that which is afforded the monstrous McKay--Thompson series (\ref{eqn:MTseries}) by virtue of the fact that they are normalized principal moduli (\ref{eqn:normalizedpm}). Since the rate of growth in the Fourier coefficients of a weakly holomorphic modular form 
is determined by its behavior 
near cusps (cf. e.g. \cite{Dabholkar:2012nd}) 
we refer to the property (\ref{eqn:Thpmanalogue}) as {\em optimality} for the McKay--Thompson series of Thompson moonshine. This definition is inspired by the notion of optimality in \cite{mum} (cf. \S\ref{sec:umbralmoonshine}, (\ref{eqn:umpmanalogue})), which was in turn motivated by the optimal growth condition formulated in \cite{Dabholkar:2012nd}.

For completeness we mention that the theta series arising in the formulation of optimality (\ref{eqn:Thpmanalogue}) for Thompson moonshine are linear combinations of the functions $\theta(d^2\tau)$, for integers $d$, where $\theta(\tau):=(\theta_{0}(\tau),\theta_{1}(\tau))$ and 
\begin{gather}\label{eqn:theta}
\theta_{r}(\tau):=\sum_{n\equiv r\xmod 2}q^{\frac{n^2}{4}}.
\end{gather}

\section{Traces of Singular Moduli}\label{sec:tracesofsingularmoduli}

We have seen in the last section that 
the $j$-invariant has both arithmetic properties and an attendant sporadic simple group, even when regarded as a function on $\HH$. However the connection to rational points on elliptic curves---i.e. elliptic curve arithmetic---is perhaps still obscure. To illuminate it we seek a more structured understanding of weighted sums of singular moduli such as (\ref{eqn:85995asatrace}). 

To setup a framework for this 
we recall that
a {\em binary quadratic form} is a polynomial 
\begin{gather}\label{eqn:QABC}
Q(x,y)=Ax^2+Bxy+Cy^2
\end{gather} 
with $A,B,C\in\ZZ$.
The {modular group} $\SL_2(\ZZ)$ 
acts naturally on binary quadratic forms via
\begin{gather}\label{eqn:Qslashgamma}
\left(Q\left|\left(\begin{matrix}a&b\\c&d\end{matrix}\right)\right.\right)(x,y):=Q(ax+by,cx+dy).
\end{gather}
The {\em discriminant} $D:=B^2-4AC$ of $Q$ as in (\ref{eqn:QABC}) is $\SL_2(\ZZ)$-invariant so it is natural to 
consider 
the $\SL_2(\ZZ)$-orbits on
\begin{gather}\label{eqn:mcQD}
\mc{Q}(D):=\left\{Q(x,y)=Ax^2+Bxy+Cy^2\mid D=B^2-4AC\right\}. 
\end{gather}

Indeed, 
if $D$ is square-free and $D\equiv 1\xmod 4$, or if $D= 4d$ for some square-free $d$ such that $d\equiv 2,3\xmod 4$, 
then $D$ is the discriminant of the number field $\QQ(\sqrt{D})$, and the orbit space $\mc{Q}(D)/\SL_2(\ZZ)$ is in natural correspondence with 
the narrow ideal class group of $\QQ(\sqrt{D})$ for $D>0$, or with $2$ copies of the ideal class group of $\QQ(\sqrt{D})$ for $D<0$ (see e.g. Chapter 5 of \cite{MR1228206} or Chapter 6 of \cite{MR1790423}).

Call a non-zero integer $D$ a {\em discriminant} if $\mc{Q}(D)$ is not empty (i.e. if $D\equiv 0,1\xmod 4$), and call $D$ a {\em fundamental discriminant} if $D$ is the discriminant of $\QQ(\sqrt{D})$. 
Define $h(D)$ to be the order of the ideal class group of $\QQ(\sqrt{D})$. Then according to the above we have  
\begin{gather}\label{eqn:hD}
	h(D)=\frac12\#\mc{Q}(D)/\SL_2(\ZZ)
\end{gather}
for $D$ negative and fundamental.

Suppose now that $f$ is a 
holomorphic function on $\HH$. Say that $f$ is a {\em weakly holomorphic modular form of weight $0$} for $\SL_2(\ZZ)$ if it is $\SL_2(\ZZ)$-invariant and satisfies $f(\tau)=O(e^{C\Im(\tau)})$ as $\Im(\tau)\to \infty$ for some $C>0$. 
Then for $D$ a negative discriminant the 
associated {\em trace of singular moduli} is
\begin{gather}\label{eqn:trDf}
	\tr(f|D):=\sum_{Q\in \mc{Q}(D)/\SL_2(\ZZ)}\frac{f(\tau_Q)}{\#\SL_2(\ZZ)_Q}
\end{gather}
where in each summand $\tau_Q$ is the unique root of $Q(x,1)=0$ with $\Im(\tau_Q)>0$ and $\SL_2(\ZZ)_Q$ is the subgroup of $\SL_2(\ZZ)$ that fixes $Q$. 

If $D$ is a negative discriminant and $D_0$ is a fundamental 
discriminant divisor of $D$ a corresponding {\em twisted trace of singular moduli} may be defined by introducing a factor $\frac{1}{\sqrt{D_0}}\chi_{D_0}(Q)$ to each summand in (\ref{eqn:trDf}), for a suitably defined function 
$\chi_{D_0}$ (see \S I.2 of \cite{MR909238} for the definition). In (\ref{eqn:85995asatrace}) we see the special case of this that $D=-15$ and $D_0=5$.

Traces of singular moduli arise as coefficients of modular forms. Indeed, one sees from the analysis in \cite{Zag_TrcSngMdl} for example that for any $f$ 
as above with vanishing constant term (i.e. $a_f(0)=0$ in the expansion $f(\tau)=\sum_n a_f(n)q^n$)
there exists a unique polynomial $\TP f(x)$ such that 
\begin{gather}\label{eqn:Fhat}
\TS f(\tau)
:=\TP f(q^{-1}) 
+\sum_{D<0}\tr(f|D)q^{|D|}
\end{gather}
is a (weakly holomorphic) modular form of weight $\frac32$ for $\Gamma_0(4)$ (\ref{eqn:Gamma0N}). In fact the same is true even if $f$ has a non-vanishing constant term, except that in that case $\TS f$ is mock modular rather than modular.

To explain what mock modular means in this context we consider the case that $f(\tau)=1$. Then $\tr(f|D)=\tr(1|D)$ is the {\em Hurwitz class number} of $D$, and we have
\begin{gather}\label{eqn:tr1DhD}
	\tr(1|D)=h(D)
\end{gather}
for $D<-4$ fundamental (since $\SL_2(\ZZ)_Q=\{\pm I\}$ for $Q\in \mc{Q}(D)$ for such $D$).
The polynomial $\TP f=\TP\,1$ is the constant polynomial $\TP\, 1(x)=-\frac1{12}$, and the function 
\begin{gather}\label{eqn:msH}
\ms{H}(\tau):=\TS\, 1(\tau)
=-\frac1{12}+\frac13q^3+\frac12q^4+q^7+q^8+\dots 
\end{gather} 
produced by (\ref{eqn:Fhat}) is the {\em Hurwitz class number generating function}, whose modular properties were first determined by Zagier \cite{MR0429750} (see also \cite{MR0453649}). In modern language 
we say that $\mathscr{H}$ is a {\em mock modular form} of weight $\frac32$ for $\Gamma_0(4)$ with {\em shadow} given by
\begin{gather}\label{eqn:checktheta}
\check\theta(\tau):=\theta_0(4\tau)+\theta_1(4\tau)=
\sum_n q^{n^2}
\end{gather} 
(cf. (\ref{eqn:theta})). In practical terms this means that the {\em completion}
\begin{gather}\label{eqn:msF}
	\widehat{\ms{H}}(\tau):=
	\ms{H}(\tau) 
	+\Im(\tau)^{-\frac12}\sum_n\beta(4\pi n^2\Im(\tau))q^{-n^2}
\end{gather}
is a real analytic modular form of weight $\frac32$ for $\Gamma_0(4)$, where 
\begin{gather}
\beta(t):=\frac{1}{16\pi}\int_1^\infty u^{-\frac32}e^{-ut}{\rm d}u.
\end{gather}

We refer the reader to \cite{Dabholkar:2012nd,Ono_unearthing} for more on mock modular forms, and to \cite{MR2248151,MR1930980} for more thorough and more general analyses of the construction (\ref{eqn:Fhat}).

\section{O'Nan Moonshine}\label{sec:ONanmoonshine}

For another example of a (non-mock) modular form obtained via traces of singular moduli consider the case that $f=f^\ON$ where
$f^\ON:=\frac12j^2-\frac{1489}2j+80256$ is the unique $\SL_2(\ZZ)$-invariant holomorphic function on $\HH$ that satisfies 
\begin{gather}\label{eqn:fON}
f^\ON(\tau)=\frac12q^{-2}-\frac12q^{-1}+O(q)
\end{gather}
as $\Im(\tau)\to\infty$. Then $\TP {f^\ON}(x)=2-x^4$ and 
\begin{gather}\label{eqn:fONhat}
 \TS f^\ON(\tau) = -q^{-4}+2+26752q^{3}+143376q^4+8288256q^7+\dots
\end{gather}
(cf. (\ref{eqn:Fhat})) is a weakly holomorphic modular form of weight $\frac32$ for $\Gamma_0(4)$.

At this point the reader may not be surprised to learn of a connection between (\ref{eqn:fONhat}) and a sporadic simple group. The group in question was discovered by O'Nan in 1973 (cf. \cite{MR0401905}) 
and has order 
\begin{gather}
\#\ON=460815505920. 
\end{gather} 
(According to the classification of finite simple groups \cite{MR2072045}, most finite simple groups, and the sporadic simple groups in particular, are uniquely determined amongst finite simple groups by their orders. See \cite{MR1023806}.) The O'Nan group is called a {\em non-monstrous} or {\em pariah} sporadic group because, according to \cite{MR671653}, it is not one of the $20$ sporadic simple groups that occur as a quotient of a subgroup of the monster.

{\em O'Nan moonshine} \cite{2017arXiv170203516D,2017NatCo...8..670D} begins with the observation that $26752$ 
is the dimension of an irreducible representation of $\ON$ (cf. \cite{atlas}). One then observes that there are analogues of $\TS f^\ON$ for certain subgroups of $\Gamma_0(4)$ that reproduce character values of non-trivial elements of $\ON$ for this representation. For example, there is a natural analogue of (\ref{eqn:fONhat}) 
for $\Gamma_0(12)=\Gamma_0(4\cdot 3)$ that satisfies $-q^{-4}+2+O(q)$ as $\Im(\tau)\to \infty$ (cf. (\ref{eqn:ONpmanalogue})). The coefficient of $q^3$ in this form is $22$, and $22$ is precisely the value of the unique irreducible $26752$-dimensional character of $\ON$ on any element of order $3$. 

Ultimately we arrive at an assignment---see Theorem 4.1 in \cite{2017arXiv170203516D}---of weakly holomorphic modular forms 
\begin{gather}\label{eqn:fONcheck}
	\check F^\ON_g(\tau)=-q^{-4}+2+\sum_{\substack{D< 0\\D\equiv 0,1\xmod 4}}C^\ON_g(D)q^{|D|}
\end{gather}
of weight $\frac32$ for $\Gamma_0(4o(g))$ to elements $g\in \ON$, with $\check F^\ON_e=\TS f^\ON$ (\ref{eqn:fONhat}), for which $g\mapsto C^\ON_g(D)$ is a {\em virtual} character of $\ON$ (i.e. an integer combination of irreducible characters of $\ON$) for each fixed $D<0$.

In this case too we find an analogue of the principal modulus property (\ref{eqn:normalizedpm}) of monstrous moonshine, which may also be compared to the corresponding optimality property (\ref{eqn:Thpmanalogue}) of Thompson moonshine. Again it is best to formulate it in terms of vector-valued forms, so we define $F^\ON_g:=(F^\ON_{g,0},F^\ON_{g,1})$ where, similar to (\ref{eqn:vvFTH}), we set
\begin{gather}\label{eqn:vvFON}
	F^\ON_{g,r}(\tau):=\delta_{r,0}(-q^{-1}+2)+\sum_{\substack{D<0\\D\equiv r^2\xmod 4}}C_g^\ON(D)q^{\frac{|D|}{4}}
\end{gather}
for $r\in \{0,1\}$.
Then for each $g\in \ON$ the corresponding vector-valued McKay--Thompson series $F_g^\ON$ is uniquely determined---up to a cusp form---amongst modular forms of weight $\frac32$ (with a suitably defined multiplier) for $\Gamma_g^{\ON}:=\Gamma_0(o(g))$ by the property that 
\begin{gather}\label{eqn:ONpmanalogue}
F_g^\ON(\tau)=
\begin{cases}
(-q^{-1},0)+O(1) &\text{ as $\Im(\tau)\to \infty$,}\\ 
O(1) &\text{ as $\tau$ tends to any non-infinite cusp of $\Gamma_g^{\ON}$.}
\end{cases}
\end{gather}

The condition (\ref{eqn:ONpmanalogue}) is similar to (\ref{eqn:Thpmanalogue}) but 
we should not refer to it simply as optimality 
since there exist modular forms of weight $\frac32$ with the same multiplier system that grow like $(0,q^{-\frac14})+O(1)$ as $\Im(\tau)\to \infty$. We refer to (\ref{eqn:ONpmanalogue}) as {\em $(-q^{-1})$-optimality} 
since it singles out the modular forms that have optimal growth subject to satisfying $(-q^{-1},0)+O(1)$ as $\Im(\tau)\to \infty$.

A {\em cusp form} is a modular form that vanishes at all cusps (\ref{eqn:cusps}). 
Cusp forms play an important role in O'Nan moonshine, and this
becomes evident when we attempt to express the McKay--Thompson series $F^\ON_g$ in terms of singular moduli. We have 
\begin{gather}
C^\ON_{e}(D)=\tr(f^\ON|D)
\end{gather}
for $D<0$ 
by definition (cf. (\ref{eqn:Fhat}), (\ref{eqn:fON}), (\ref{eqn:fONhat}), (\ref{eqn:fONcheck})). 
To obtain analogues of this 
for non-trivial elements of the O'Nan group we define
\begin{gather}\label{eqn:trNDF}
	\tr_N(f|D):=\sum_{Q\in \mc{Q}_N(D)/\Gamma_0(N)}\frac{f(\tau_Q)}{\#\Gamma_0(N)_Q},
\end{gather}
for $f$ any weakly holomorphic modular form of weight $0$ for $\Gamma_0(N)$ (\ref{eqn:Gamma0N}), where 
\begin{gather}\label{eqn:mcQND}
\mc{Q}_N(D):=\left\{Q\in \mc{Q}(D)\mid A\equiv 0\xmod N\right\}. 
\end{gather}
Then for $o(g)=3$, for example, Appendix D and Proposition 5.1 of \cite{2017arXiv170203516D} tell us that 
\begin{gather}\label{eqn:cONgo3}
C_{3A}^\ON(D)=12\tr_1(1|D)-12\tr_3(1|D)+\tr_3(f^\ON_{3A}|D)
\end{gather}
for $D<0$. Here $C^\ON_{3A}(D)$ is $C^\ON_g(D)$ for any $g\in \ON$ with $o(g)=3$ (since $\ON$ has a unique conjguacy class of elements of order $3$, denoted $3A$ in \cite{atlas}), and $f^\ON_{3A}$ is the unique $\Gamma_0(3)$-invariant holomorphic function on $\HH$ that satisfies 
\begin{gather}\label{eqn:fON3A}
f^\ON_{3A}(\tau)=\frac12q^{-2}-\frac12q^{-1}+O(q)
\end{gather}
as $\Im(\tau)\to \infty$, and remains bounded as $\tau\to 0$. 

No cusp forms can appear in (\ref{eqn:cONgo3}) because there are no non-zero cusp forms of weight $\frac32$ with the correct multiplier system (i.e., the inverse of that of $\theta(\tau)=(\theta_0(\tau),\theta_1(\tau))$ as in (\ref{eqn:theta})) for $\Gamma_0(3)$. However, there is a unique up to scale such cusp form $\ms{G}_{15}(\tau)=(\ms{G}_{15,0}(\tau),\ms{G}_{15,1}(\tau))$ for $\Gamma_0(15)$, which for a suitable choice of scaling satisfies
\begin{gather}\label{eqn:checkmsG15}
\check{\ms{G}}_{15}(\tau)=\sum_{D<0}C_{15}(D)q^{|D|}=q^3-2q^8-q^{15}+2q^{20}+\dots 
\end{gather}
where $\check{\ms{G}}_{15}(\tau):=\ms{G}_{15,0}(4\tau)+\ms{G}_{15,1}(4\tau)$
(cf. (\ref{eqn:checktheta})).
With this definition we see from Appendix D and Proposition 5.3 of \cite{2017arXiv170203516D} that
\begin{gather}\label{eqn:cONgo15}
C_{15AB}^\ON(D)=\ldots+\frac14\tr_{15}(f^\ON_{15AB}|D)+\frac94{C}_{15}(D)
\end{gather}
for $D<0$, where $C^\ON_{15AB}(D)$ is $C^\ON_g(D)$ for any $g\in \ON$ with $o(g)=15$, and where $f^\ON_{15AB}$ is now an analogue of $f^\ON$ (\ref{eqn:fON}) for $\Gamma_0(15)$. 
The ellipsis in (\ref{eqn:cONgo15}) stands in for a linear combination of traces $\tr_N(1|D)$ for $N$ ranging over the divisors of $15$.

So in particular the cusp form $\ms{G}_{15}$ makes rational, generally non-integral contributions to the coefficients of $F^\ON_g$ for $o(g)=15$. So the cusp form coefficients $C_{15}(D)$ must satisfy congruences with the particular linear combinations of traces of singular moduli that appear in (\ref{eqn:cONgo15}), since the McKay--Thompson series coefficients $C^\ON_g(D)$ are all algebraic integers by force of being (virtual) character values.
As we will see presently, it is this simple circumstance that opens the door to connecting the O'Nan group with elliptic curve arithmetic.

\section{Cusp Forms and Elliptic Curves}\label{sec:cuspformsellipticcurves}

We have alluded to the fact that cusp forms of weight $\frac32$, such as $\ms{G}_{15}$ (\ref{eqn:checkmsG15}), are connected to elliptic curve arithmetic. 
This is so because Waldspurger proved 
\cite{MR577010,MR646366} that the coefficients of such cusp forms 
may be expressed in terms of special values of $L$-functions of elliptic curves defined over $\QQ$. 
Such special values, in turn, 
constrain the groups of rational points on the corresponding elliptic curves. 

To explain this we 
consider
the 
{\em Hasse--Weil $L$-function} 
of an elliptic curve $E$ over $\QQ$, 
which is defined, for $s\in \CC$ with $\Re(s)$ sufficiently large, by 
\begin{gather}\label{eqn:LEs}
	L_E(s):=\prod_{p\text{ prime}}\frac1{1-a_{E}(p)p^{-s}+\epsilon(p)p^{1-2s}}.
\end{gather}
Here $a_{E}(p):=p+1-\#E(\FF_p)$, for $E(\FF_p)$ the group of $\FF_p$-rational points on the mod $p$ reduction of $E$, and $\epsilon(p)$ is $0$ or $1$ as $p$ divides the minimal discriminant of $E$ or not. By {\em minimal discriminant} we mean the discriminant of a {\em global minimal Weierstrass model} for $E$, 
which is an equation of the form
\begin{gather}\label{eqn:minimalWeierstrassmodel}
	y^2+a_1xy+a_3y=x^3+a_2x^2+a_4x+a_6,
\end{gather}
with $a_1,a_2,a_3,a_4,a_6\in\ZZ$, 
which (after passage to the corresponding homogeneous cubic) defines a curve in $\PP^2(\QQ)$ isomorphic to $E$ and has minimal $|\Delta|$ amongst all such equations, where $\Delta:=-b_2^2b_8-8b_4^3-27b_6^2+9b_2b_4b_6$ for 
\begin{gather}
\begin{split}
	b_2&:=a_1^2+4a_2,\\
	b_4&:=a_1a_3+2a_4,\\
	b_6&:=a_3^2+4a_6,\\
	b_8&:=a_1^2a_6-a_1a_3a_4+4a_2a_6+a_2a_3^2-a_4^2.
\end{split}
\end{gather}

It turns out (see e.g. \S3.1 of \cite{MR1628193}) that any elliptic curve $E$ defined over $\QQ$ has a unique {\em reduced} global minimal Weierstrass model (\ref{eqn:minimalWeierstrassmodel}) satisfying $a_1,a_3\in \{0,1\}$ and $a_2\in \{-1,0,1\}$. 
For example, for
\begin{gather}\label{eqn:E0}
	E_{15}\,:\, 
	y^2=x^3-12987x-263466	
\end{gather}
the reduced global minimal model is 
\begin{gather}\label{eqn:E0redmin}
	y^2+xy+y=x^3+x^2-10x-10,
\end{gather}
which we may obtain from (\ref{eqn:E0}) by replacing $x$ with $36x+15$, replacing $y$ with $108x+216y+108$, and dividing the result by $46656$.
So the minimal discriminant of $E_{15}$ is $\Delta=50625=3^4\cdot5^4$, and so $\epsilon(p)=0$ in (\ref{eqn:LEs}) just for $p=3$ and $p=5$, for $E=E_{15}$.

Note also that a global minimal Weierstrass equation for $E$ should be used when computing the mod $p$ reductions of $E$ for $p=2$ and $p=3$. 
So for example from (\ref{eqn:E0redmin}) we see that the reduction of $E_{15}$ mod $2$ is given by $y^2+xy+y=x^3+x^2$ (rather than $y^2=x^3+x$)
and the reduction of $E_{15}$ mod $3$ is given by $y^2+xy+y=x^3+x^2-x-1$ (rather than $y^2=x^3$). From this we see that
$\#E_{15}(\FF_2)=4$ and $\#E_{15}(\FF_3)=5$, 
so we have
\begin{gather}\label{eqn:a2a3E0}
a_{E_{15}}(2)=a_{E_{15}}(3)=-1.
\end{gather}

The significance of 
$L_E$ (\ref{eqn:LEs})
for the computation of $E(\QQ)$ 
is explicated by the {\em Birch--Swinnerton-Dyer conjecture} \cite{MR0179168}, 
which predicts that the rank of $E(\QQ)$ (as a $\ZZ$-module) is precisely the order of vanishing of $L_E(s)$ at $s=1$. This is the {\em weak form} of the conjecture. The {\em strong form} states that if $E(\QQ)\simeq \ZZ^r\oplus E(\QQ)_\tor$, where $E(\QQ)_\tor$ is the subgroup of finite order elements of $E(\QQ)$, then 
\begin{gather}\label{eqn:strongBSD}
	\lim_{s\to 1}\frac{L^{(r)}_E(s)}{r!\Omega_E } = C_E\frac{\#\Sha(E)}{(\#E(\QQ)_\tor)^2}
\end{gather}
where $L^{(n)}_E(s)$ denotes the $n$-th derivative of $L_E(s)$, we write $\Omega_E$ and $C_E$ for certain computable constants attached to $E$, and $\Sha(E)$ is the {\em Tate--Shafarevich group} of $E$. 

We refer to 
\cite{MR2514094} for more detail on 
the invariants we have denoted $\Omega_E$, $C_E$ and $\Sha(E)$. The most subtle of these is the Tate--Shafarevich group $\Sha(E)$, which is not yet known to be finite except in special cases. So in particular the right-hand side of (\ref{eqn:strongBSD}) is not generally known to be well-defined. Roughly speaking, $\Sha(E)$ encodes the obstruction to recovering $E(\QQ)$ from computations with $E$ modulo primes (cf. e.g. \S7 of \cite{MR0419359}). 
An important related object is the {\em $\ell$-Selmer group} of $E$, 
for $\ell$ an integer, 
which satisfies a short exact sequence
\begin{gather}\label{eqn:SelmerSES}
	0\to E(\QQ)/\ell E(\QQ)\to \Sel_\ell(E)\to \Sha(E)[\ell]\to 0.
\end{gather}
Here $\Sha(E)[\ell]$ denotes the elements of order dividing $\ell$ in $\Sha(E)$. 
The $\ell$-Selmer group of an elliptic curve over $\QQ$ is known to be finite for every $\ell$ (cf. e.g. \cite{MR2514094}), so at least the $\ell$-torsion in $\Sha(E)$ is finite for every $\ell$.

Interestingly there is subtlety involved in defining the left-hand side of (\ref{eqn:strongBSD}) also, for until the turn of the century the best known bound on the $a_E(p)$ in (\ref{eqn:LEs}) was Hasse's 1933 result (established, amongst other things, in \cite{MR1581496,MR1581499,MR1581508}) that 
\begin{gather}
|a_E(p)|\leq 2\sqrt{p}.
\end{gather} 
This is enough to prove that the product in (\ref{eqn:LEs}) converges for $\Re(s)>\frac32$, but not enough to define $L_E(s)$ near $s=1$. 

The fact that $L_E$ can be continued to the entire complex plane follows from the {\em modularity theorem}. 
This result implies that for any elliptic curve $E$ over $\QQ$ there exists a cusp form $f$ of weight $2$ for $\Gamma_0(N)$, for some positive integer $N=N_E$, satisfying 
\begin{gather}\label{eqn:FEFricke}
f\left(-\frac1{N\tau}\right)=w_NN\tau^2f(\tau)
\end{gather} 
for some $w_N\in \{\pm 1\}$, 
such that $L_E(s)=L_f(s)$ where $\frac{\Gamma(s)}{(2\pi)^{s}}L_f(s)$ is the 
{\em Mellin transform} of $f(it)$, 
\begin{gather}\label{eqn:Lfs}
	L_f(s) := \frac{(2\pi)^s}{\Gamma(s)}\int_0^\infty f(it) t^{s-1}{{\rm d}t}.
\end{gather}
Here $\Gamma(s)$ denotes the {\em Gamma function} $\Gamma(s):=\int_0^\infty e^{-t}t^{s-1}{{\rm d}t}$ (which is the Mellin transform of $e^{-t}$). So in particular 
we should have
\begin{gather}\label{eqn:aEpafp}
a_E(p)=a_f(p)
\end{gather}  
for $p$ prime (cf. (\ref{eqn:LEs})), when $f(\tau)=\sum_{n}a_f(n)q^n$.
Given such a cusp form $f$ it is (\ref{eqn:FEFricke}) that allows one to continue $L_E$ to the entire complex plane. This is because (\ref{eqn:FEFricke}) implies, via (\ref{eqn:Lfs}),  a corresponding {\em functional equation} 
that relates $L_E(s)$ to $L_E(2-s)$.

The full proof of the modularity theorem was announced in 1999 by Breuil--Conrad--Diamond--Taylor (cf. \cite{MR1723249}). Their work \cite{MR1839918} built upon earlier work of Diamond \cite{MR1405946} and Conrad--Diamond--Taylor \cite{MR1639612}, which in turn built upon the breakthrough results of Wiles \cite{MR1333035} and Taylor--Wiles \cite{MR1333036} that famously verified Fermat's last theorem. (The work of Wiles on Fermat's last theorem is beautifully exposited in \cite{MR3699470,MR3617390}.)

For an example of the modularity theorem in action consider $E=E_{15}$ as in (\ref{eqn:E0}). In this case $N=15$ and the corresponding cusp form $f=f_{15}$ 
happens to admit a product formula (this is not typical),
\begin{gather}\label{eqn:FE0}
	f_{15}(\tau)
	=\eta(\tau)\eta(3\tau)\eta(5\tau)\eta(15\tau)
	=q-q^2-q^3-q^4+\dots.
\end{gather}
Here $\eta(\tau)$ is as in (\ref{eqn:jEtau}).
According to the {\em Euler identity} we have
\begin{gather}\label{eqn:EulerId}
	\eta(\tau)=\sum_{n>0} \left(\tfrac{12}n\right)q^{\frac{n^2}{24}}
\end{gather}
where $\left(\tfrac{12}n\right)$ is $1$ when $n\equiv 1,11\xmod 12$, and $-1$ when $n\equiv 5,7\xmod 12$, and zero otherwise. The series (\ref{eqn:EulerId}) converges fast enough to make it effective for numerical computation.

We have $f_{{15}}(\gamma\tau){\rm d}\gamma\tau=f_{{15}}(\tau){\rm d}\tau$ for $\gamma\in \Gamma_0(15)$, which, together with vanishing at cusps, is just what it means for $f_{{15}}(\tau)$ to be a cusp form of weight $2$ for $\Gamma_0(15)$. 
Comparing (\ref{eqn:a2a3E0}) to (\ref{eqn:FE0}) we see that (\ref{eqn:aEpafp}) is satisfied at least for $p=2,3$. One can check that $w_{15}=-1$ in (\ref{eqn:FEFricke}), and this implies that the order of vanishing of $L_{E_{15}}(s)$ at $s=1$ is even. So if the Birch--Swinnerton-Dyer conjecture is true for $E_{15}$ then $E_{15}(\QQ)$ has even rank. It can be shown that in fact 
\begin{gather}\label{eqn:E0(Q)}
E_{15}(\QQ)\simeq \ZZ/2\ZZ\times \ZZ/4\ZZ
\end{gather} 
has rank $0$ (cf. e.g. \cite{lmfdb}). Using (\ref{eqn:Lfs}), (\ref{eqn:FE0}) and (\ref{eqn:EulerId}) one can compute numerically that $L_{E_{15}}(1)\approx 0.3501507606$.

There are a number of results in the direction of (\ref{eqn:strongBSD}) but the Birch--Swinnerton-Dyer conjecture is largely open (in both the strong and weak forms, cf. \cite{MR2238272}). It is known that if $L_E(1)\neq 0$ then $E(\QQ)$ is finite, and if $L_E(1)=0$ but $L'_E(1)\neq 0$ then $E(\QQ)$ has rank $1$. The former statement follows from the modularity theorem together with work of Kolyvagin \cite{MR954295}. The later statement follows from modularity, Kolyvagin's work loc. cit., and the construction by Gross--Zagier \cite{MR833192} of a point of infinite order in $E(\QQ)$ in case $L_E(s)$ vanishes to order $1$ at $s=1$. Kolyvagin's work also shows that $\Sha(E)$ is finite in these cases. Significant for the sequel is a $p$-local version of (\ref{eqn:strongBSD}) that was obtained by Skinner \cite{MR3513846} following earlier work \cite{MR3148103} of Skinner--Urban.

\section{Elliptic Curves and O'Nan}\label{sec:ECONan}

We are now ready to present a 
result that 
relates the O'Nan sporadic simple group to elliptic curves over the rationals. For this let $E_{15}\otimes D$, for $D$ a non-zero integer, be the elliptic curve over $\QQ$ defined by
\begin{gather}\label{eqn:E15D}
	E_{15}\otimes D \,:\, y^2 = x^3 - 12987D^2x -  263466D^3.
\end{gather}
This is the {\em $D$-th quadratic twist} of $E_{15}=E_{15}\otimes 1$ (\ref{eqn:E0}). In general $E_{15}\otimes D$ has the same $j$-invariant as $E_{15}$ and is isomorphic to $E_{15}$ over $\QQ(\sqrt{D})$. There is a notion of Hasse--Weil $L$-function (\ref{eqn:LEs}) for elliptic curves over an arbitrary number field, and a corresponding formulation of the Birch--Swinnerton-Dyer conjecture (\ref{eqn:strongBSD}). Taking $D$ to be fundamental 
(cf. \S\ref{sec:tracesofsingularmoduli})
one finds that the $L$-function for $E$ over $\QQ(\sqrt{D})$ is the product of the $L$-functions for $E$ and $E\otimes D$, each regarded as curves over $\QQ$. So, conjecturally at least, the behavior of $L_{E_{15}\otimes D}(s)$ near $s=1$ also controls the arithmetic of $E_{15}$ 
as a curve over $\QQ(\sqrt{D})$.

The significance of (\ref{eqn:E15D}) for us is that if $D$ is a negative fundamental discriminant then the aforementioned result of 
Waldspurger (cf. the first paragraph of \S\ref{sec:cuspformsellipticcurves}) relates the $L$-function special value $L_{E_{15}\otimes D}(1)$ (\ref{eqn:LEs}) 
to the cusp form coefficient $C_{15}(D)$ (\ref{eqn:checkmsG15}) 
that 
appears in the expression 
(\ref{eqn:cONgo15}) 
for the 
McKay--Thompson series coefficient $C^\ON_{15AB}(D)$ (\ref{eqn:fONcheck}) of moonshine for O'Nan.

In joint work \cite{2017arXiv170203516D} with Michael Mertens and Ken Ono we use an explicit formulation of Waldspurger's result due to Kohnen \cite{MR783554}, together with a $p$-local version of (\ref{eqn:strongBSD}) due to Skinner \cite{MR3513846} and Skinner--Urban \cite{MR3148103}, to establish the following theorem, which intertwines the character values $C^\ON_{3A}(D)$ (\ref{eqn:cONgo3}), the class numbers $h(D)$ (\ref{eqn:hD}), and the $5$-Selmer and Tate--Shafarevich groups of the elliptic curves $E_{15}\otimes D$ (\ref{eqn:E15D}). 
\begin{theorem}[\cite{2017arXiv170203516D}]
\label{thm:ONEC}
Let $D$ be a negative fundamental discriminant such that $D\equiv 1\xmod 3$ and $D\equiv 2,3\xmod 5$. 
Then $\Sel_5(E_{15}\otimes D)\neq\{0\}$ if and only if 
\begin{gather}\label{eqn:ONECcong}
C_{3A}^\ON(D)+h(D)\equiv 0\xmod 5. 
\end{gather}
Furthermore, if $D$ is as above and $L_{E_{15}\otimes D}(1)\neq 0$ then 
$\#\Sha(E_{15}\otimes D)\equiv 0\xmod 5$ if and only if 
(\ref{eqn:ONECcong}) holds.
\end{theorem}

The short exact sequence (\ref{eqn:SelmerSES}) tell us that if $\Sel_\ell(E)$ is trivial for any $\ell>1$ then $E(\QQ)$ is finite. 
So Theorem \ref{thm:ONEC} gives a new criterion for checking when a quadratic twist $E_{15}\otimes D$ (\ref{eqn:E15D}) of $E_{15}$ (\ref{eqn:E0}) has finitely many rational points.

Another feature of Theorem \ref{thm:ONEC} is that the 
criterion (\ref{eqn:ONECcong}) does not 
involve cusp forms, or any of the $E_{15}\otimes D$. It can be checked for any admissible $D$ as soon as we know a principal modulus for $\Gamma_0(3)$ (cf. (\ref{eqn:cONgo3}), (\ref{eqn:fON3A})), and enough about $\mc{Q}(D)$ (\ref{eqn:mcQD}). 
So in a sense, Theorem \ref{thm:ONEC} says that the 
O'Nan group ``reduces'' the cubic equation (\ref{eqn:E15D}) defining $E_{15}\otimes D$ to the quadratic equations $Ax^2+Bx+C=0$ with $B^2-4AC=D$.

Let us put Theorem \ref{thm:ONEC} into practice for $D=-8$ and $D=-68$. For these values of $D$ the corresponding elliptic curves are 
\begin{gather}
E_{15}\otimes (-8)=E_{15}',\quad
E_{15}\otimes(-68)=E_{15}'',
\end{gather}
which appeared already (\ref{eqn:E1E2}) in \S\ref{sec:ellipticcurves}. Using Theorem \ref{thm:ONEC} we will see $\Sel_5(E_{15}')$ is trivial but $\Sel_5(E_{15}'')$ is not. So in particular we will confirm the claim of \S\ref{sec:ellipticcurves} that $E_{15}'(\QQ)$ is finite.

We first compute $h(-8)$ and $h(-68)$. For $D<-4$ fundamental we have $h(D)=\tr(1|D)$ (cf. (\ref{eqn:tr1DhD})), and we can compute $\tr(1|D)$ 
directly 
by considering {reduced} quadratic forms. The idea here is to constrain the coefficients $A,B,C$ in (\ref{eqn:QABC}) so that $\tau_Q=\frac{-B+\sqrt{D}}{2A}$ lies in a fundamental domain for $\SL_2(\ZZ)$. Concretely, say that $Q(x,y)=Ax^2+Bxy+Cy^2$ 
with $D<0$
is {\em reduced} if 
\begin{gather}\label{eqn:reducedABC}
|B|\leq A\leq C,
\end{gather}
and if $B\geq 0$ when $A=|B|$ or $A=C$. Then there is exactly one reduced quadratic form in each orbit of $\SL_2(\ZZ)$ on {\em positive-definite} forms in $\mc{Q}(D)$ (i.e. those with $A>0$). 
From (\ref{eqn:reducedABC}) we have $|D|=4AC-B^2\geq 4A^2-A^2=3A^2$. So 
the problem of computing $h(D)$ reduces, via reduced forms, to consideration of the pairs $(A,B)$ such that $3A^2\leq {{|D|}}$ and $-A<B\leq A$ and $B^2\equiv D\xmod 4A$.
Given such a pair $(A,B)$ we obtain a reduced quadratic form of discriminant $D$ by setting $C=\frac{B^2-D}{4A}$. 

We have $0<3A^2\leq 8$ only for $A=1$. For $-1< B\leq 1$ we have $B^2\equiv -8\xmod 4$ only for $B=0$, so there is exactly one reduced form with $D=-8$ (in agreement with (\ref{eqn:msH})). 
It only takes a little more work to see that for $D=-68$ the relevant pairs are $(1,0)$, $(2,2)$ and $(3,\pm 2)$, so 
\begin{gather}\label{eqn:h8h68}
	h(-8)=1,\quad h(-68)=4.
\end{gather}
(For a comprehensive discussion of class number computation see Chapter 5 of \cite{MR1228206}.)

To compute $C^\ON_{3A}(D)$ via traces (\ref{eqn:cONgo3}) we require representatives for the orbits of $\Gamma_0(3)$ on $\mc{Q}_3(D)$ (\ref{eqn:mcQND}). For this define 
\begin{gather}\label{eqn:mcQNDr}
\mc{Q}_N(D,r):=\left\{Q\in \mc{Q}_N(D)\mid B\equiv r\xmod 2N\right\}.
\end{gather}
Then $\mc{Q}_N(D,r)$ is a $\Gamma_0(N)$-invariant subset of $\mc{Q}_N(D)$, and for $(D,r)$ such that $D<0$ is fundamental and $r^2\equiv D\xmod 4N$ the Proposition in \S I.1 of \cite{MR909238} tells us that both the inclusion map $\mc{Q}_N(D,r)\hookrightarrow\mc{Q}(D)$, and the map 
\begin{gather}\label{eqn:mcQNDrmap}
Ax^2+Bxy+Cy^2\mapsto \frac{A}Nx^2+Bxy+NCy^2,
\end{gather}
induce bijections from 
$\mc{Q}_N(D,r)/\Gamma_0(N)$ to $\mc{Q}(D)/\SL_2(\ZZ)$. 

Note that $\mc{Q}_N(D,r)$ only depends on $r\xmod 2N$. For $N=3$ and $D\equiv 4\xmod 12$ the relevant $r$ for (\ref{eqn:mcQNDrmap}) are $r\equiv\pm 2\xmod 6$, so for such $D$ we have $\mc{Q}_3(D)=\mc{Q}_3(D,-2)\cup\mc{Q}_3(D,2)$. So 
$\tr_3(1|D)=2\tr_1(1|D)=2h(D)$
for $D<0$ fundamental such that $D\equiv 4\xmod 12$. (We have $\Gamma_0(N)_Q=\{\pm I\}$ for all $N$ when $Q\in \mc{Q}(D)$ and $D<-4$.) Thus (\ref{eqn:cONgo3}) reduces to 
\begin{gather}\label{eqn:cON3AD868}
C^\ON_{3A}(D)=-12h(D)+\tr_3(f^\ON_{3A}|D)
\end{gather} 
for $D=-8$ and $D=-68$.

We have seen already (\ref{eqn:h8h68}) that $\SL_2(\ZZ)$ acts transitively on the positive-definite forms in $\mc{Q}(-8)$, 
so $\Gamma_0(3)$ must act transitively on the positive-definite forms in $\mc{Q}_3(-8,-2)$ and $\mc{Q}_3(-8,2)$. We have $3x^2\pm 2xy+y^2\in \mc{Q}_3(-8,\pm 2)$,  
so from (\ref{eqn:cON3AD868}) we obtain
\begin{gather}\label{eqn:cON3A8}
	C^\ON_{3A}(-8)=-12 + \frac12 f^\ON_{3A}\left(\frac{2+\sqrt{-8}}{6}\right)
	+ \frac12 f^\ON_{3A}\left(\frac{-2+\sqrt{-8}}{6}\right).
\end{gather}

For the $D=-68$ counterpart to (\ref{eqn:cON3A8}) we apply modular transformations (\ref{eqn:Qslashgamma}) to the forms we found in the course of computing $h(-68)$  (cf. (\ref{eqn:h8h68})). For example, for $(A,B)=(1,0)$ the corresponding quadratic form is $x^2+17y^2$, which is mapped to $x^2\pm2xy+18y^2$ when we take 
$\left(\begin{smallmatrix}a&b\\c&d\end{smallmatrix}\right)=\left(\begin{smallmatrix}1&\pm1\\0&1\end{smallmatrix}\right)$ in (\ref{eqn:Qslashgamma}). These forms, in turn,  pull back to $3x^2\pm2xy+6y^2\in \mc{Q}_3(-68,\pm 2)$ under (\ref{eqn:mcQNDrmap}). Thus we obtain the $d=1$ terms in 
\begin{gather}\label{eqn:cON3A68}
	C^\ON_{3A}(-8)=-48 + \sum_{d|6}
	\left(\frac12 f^\ON_{3A}\left(\frac{2+\sqrt{-68}}{6d}\right)
	+ \frac12 f^\ON_{3A}\left(\frac{-2+\sqrt{-68}}{6d}\right)\right),
\end{gather}
and the rest may be obtained in a similar way. (See \S A of \cite{Harvey:2015mca} for some similar calculations.)

It remains to compute the values $f^\ON_{3A}(\tau_Q)$ explicitly, for the $\tau_Q$ that appear in (\ref{eqn:cON3A8}) and (\ref{eqn:cON3A68}). We have $f^\ON_{3A}(\tau)=\frac12T_3(\tau)^2-\frac12T_3(\tau)-54$ (cf. (\ref{eqn:fON3A}), (\ref{eqn:T3})) where $T_3$ is the normalized principal modulus (\ref{eqn:normalizedpm}) for $\Gamma_0(3)$. 
We also have 
\begin{gather}\label{eqn:T3}
T_3(\tau)=\frac{\eta(\tau)^{12}}{\eta(3\tau)^{12}}+12
=q^{-1}+54q-76q^2-243q^3+\dots
\end{gather}
where 
$\eta(\tau)$ is as in (\ref{eqn:jEtau}), (\ref{eqn:EulerId}). Since we know that the $C^\ON_{3A}(D)$ are integers a suitable truncation of the series (\ref{eqn:EulerId}) may be used to perform the computation.

Alternatively, we can compute the 
$C^\ON_{3A}(D)$ using $(-q^{-1})$-optimality (\ref{eqn:ONpmanalogue}) since the McKay--Thompson series $F^\ON_{g}$ 
are uniquely determined by this condition when $o(g)=3$. 
Indeed, with some work we find that
\begin{gather}\label{eqn:FON3A}
\begin{split}
F^\ON_{3A,0}(\tau)\theta_0(\tau)&+F^\ON_{3A,1}(\tau)\theta_1(\tau)=q\frac{{\rm d}}{{\rm d}q}T_3(\tau),\\
\frac{F^\ON_{3A,0}(\tau)}{\theta_0(\tau)^3}&+\frac{F^\ON_{3A,1}(\tau)}{\theta_1(\tau)^3}=
-\frac{(T_6(\tau)-5)^2(T_6(\tau)+7)(T_6(\tau)^2+\tfrac{25}4T_6(\tau)+\tfrac{35}4)}{(T_6(\tau)+3)^2(T_6(\tau)+4)^2},
\end{split}
\end{gather}
where 
$T_3$ is as in (\ref{eqn:T3}) 
and $T_6(\tau)=\frac{\eta(\tau)^5\eta(3\tau)}{\eta(2\tau)\eta(6\tau)^5}+5$ is the normalized principal modulus for 
$\Gamma_0(6)$. 
Solving 
for $F^\ON_{3A,0}(\tau)$ in (\ref{eqn:FON3A}) 
we obtain
\begin{gather}
	F^\ON_{3A,0}(\tau)
	=-q^{-1}+2+6q-188q^2+\cdots-15834144q^{17}+O(q^{18}).
\end{gather}
So we have $C^\ON_{3A}(-8)=-188$ and $C^\ON_{3A}(-68)=-15834144$, so from (\ref{eqn:h8h68}) we obtain
\begin{gather}
\begin{split}
C^\ON_{3A}(-8)&+h(-8)=-187,\\
C^\ON_{3A}(-68)&+h(-68)=-15834140.
\end{split}
\end{gather}

Thus the congruence condition (\ref{eqn:ONECcong}) of Theorem \ref{thm:ONEC} is not satisfied for $D=-8$, but is satisfied for $D=-68$, so $\Sel_5(E_{15}')$ is trivial but $\Sel_5(E_{15}'')$ is not, as we claimed. In particular, $E_{15}'(\QQ)$ is finite.

The non-triviality of $\Sel_5(E_{15}'')$ may reflect the infinitude of $E_{15}''(\QQ)$, or may be due to the existence of elements of order $5$ in $\Sha(E_{15}'')$. There are methods (other than Theorem \ref{thm:ONEC}) that confirm that the former is the case. For the sake of completeness we remark that the rank of $E_{15}''(\QQ)$ is two (cf. \cite{lmfdb}), and some independent infinite order generators are given by $(852,179712)$ and $(-3468,499392)$.

We should mention that O'Nan moonshine knows more about elliptic curve arithmetic than that which we have explained herein. How much it knows is (as far as we know) constrained by the specific numbers $N$ 
which occur as orders of elements of the O'Nan group. Here we have focused on $N=15$, 
but there are also elements of order $N=14$ to which a similar analysis may be applied. Thus in \cite{2017arXiv170203516D} we also obtain a counterpart to Theorem \ref{thm:ONEC} that relates the coefficients of the McKay--Thompson series for 
elements of order $2$ in the O'Nan group to the groups $\Sel_7(E_{14}\otimes D)$ and $\Sha(E_{14}\otimes D)[7]$ where 
\begin{gather}\label{eqn:E14D}
	E_{14}\otimes D\,:\, y^2 = x^3 + 5805D^2x -  285714D^3.
\end{gather}

\section{Umbral Moonshine}\label{sec:umbralmoonshine}

We conclude this work with an explanation of how umbral moonshine fits into the framework described above.

It was a key point in \S\ref{sec:thompsonmoonshine} that the coefficients of $q^{\frac54}$ in the McKay--Thompson series (\ref{eqn:vvFTH}) of Thompson moonshine realize the character of the direct sum of the two $85995$-dimensional irreducible representations of $\Th$, and that these characters take values in $\QQ(\sqrt{-15})$ (cf. (\ref{eqn:jtau85995})). This admits a more general formulation. Indeed, in \S4 of \cite{Harvey:2015mca} it appears as a special case of the {\em discriminant property} of Thompson moonshine, which states (amongst other things) that directly similar statements hold for the coefficients of $q^{\frac84}$ and $q^{\frac{13}4}$, and the number fields $\QQ(\sqrt{-24})$ and $\QQ(\sqrt{-39})$, respectively. 

As explained in loc. cit., an analogous discriminant property was formulated earlier for Mathieu moonshine, and for umbral moonshine more generally in \cite{um,mum}. We recall here that {\em umbral moonshine} \cite{um,mum,MR3766220,umrec} is the assignment of a vector-valued (weakly holomorphic) mock modular form $H^{(\ell)}_g(\tau)=(H^{(\ell)}_{g,r}(\tau))$,
\begin{gather}\label{eqn:vvHUM}
	H^{(\ell)}_{g,r}(\tau):=-\delta_{r,1}2q^{-\frac1{4m}}+\sum_{\substack{D\leq 0\\D\equiv r^2\xmod 4}} C^{(\ell)}_{g,r}(D)q^{\frac{|D|}{4m}},
\end{gather}
of weight $\frac12$ for $\Gamma^{(\ell)}_g:=\Gamma_0(o(g))$ 
to each 
$g\in G^{(\ell)}$, 
for a certain finite group $G^{(\ell)}$, for each $\ell$ in the following list of symbols, 
\begin{gather}\label{eqn:lambencies}
\begin{split}
	2,\; 3,\; 4,\; 5,\; 6,\; 7,\; 8,\; 9,&\; 10,\; 12,\; 13,\; 16,\; 18,\; 25, \\
	6+3,\; 10+5,\; 14+7,\; 18+&9,\; 22+11,\; 30+15,\; 46+23,\; \\
	12+4,\; 30&+6, 10,15.
\end{split}	
\end{gather}
(See \S2 of \cite{MR3766220} or \S9 of \cite{mnstmlts} for the index sets of the $r$ in (\ref{eqn:vvHUM}).) 

These 23 symbols (\ref{eqn:lambencies}) 
are called the {\em lambencies} of umbral moonshine. 
To each {lambency} 
$\ell$ 
is attached an even unimodular positive-definite lattice $N^{(\ell)}$ of rank $24$---these were classified by Niemeier \cite{Nie_DefQdtFrm24} in 1973---and the corresponding {\em umbral group} is 
\begin{gather}\label{eqn:Gell}
G^{(\ell)}:=\Aut(N^{(\ell)})/W^{(\ell)}.
\end{gather} 
Here $W^{(\ell)}$ is the normal subgroup of $\Aut(N^{(\ell)})$ generated by reflections in {\em roots} (i.e. vectors $\alpha \in N^{(\ell)}$ such that $(\alpha,\alpha)=2$). 
The shadows of the $H^{(\ell)}_g$ are (vector-valued) linear combinations of the weight $\frac32$ theta series
\begin{gather}\label{eqn:theta1mr}
\theta^1_{m,r}(\tau):=\sum_{n\equiv r\xmod 2m}n q^{\frac{n^2}{4m}}
\end{gather}
where $m=m^{(\ell)}$ 
is the positive integer for which $\ell=m+e,f,\dots$. 
The case that $\ell=2$ recovers the
{\em Mathieu moonshine} 
that was first observed by Eguchi--Ooguri--Tachikawa \cite{Eguchi2010} in 2010. For a more detailed review of umbral moonshine see \cite{mnstmlts}.

Similar to the situation in Thompson moonshine, each umbral mock modular form $H^{(\ell)}_g$ is uniquely determined---up to a theta series---amongst mock modular forms of weight $\frac12$ (with its multiplier system) by the condition that 
\begin{gather}\label{eqn:umpmanalogue}
	H^{(\ell)}_{g}(\tau)=
	\begin{cases}
		(-2q^{-\frac1{4m}},0,\ldots)+O(1)&\text{ as $\Im(\tau)\to \infty$},\\
		O(1)&\text{ as $\tau$ tends to any non-infinite cusp of $\Gamma_g^{(\ell)}$.}
	\end{cases}
\end{gather}

This 
is called {\em optimality} for umbral moonshine. 
By now the reader will readily recognize 
(\ref{eqn:umpmanalogue}) as a natural counterpart to (\ref{eqn:normalizedpm}), (\ref{eqn:Thpmanalogue}) and (\ref{eqn:ONpmanalogue}). 
(Curiously, 
optimality actually determines the $H^{(\ell)}_g$ 
uniquely except when $\ell=9$ and $o(g)\equiv 0\xmod 3$, simply because the relevant spaces of theta series turn out to be vanishing in all the other cases. This is proven in \cite{MR3766220}.)

According to \S6.4 of \cite{mum} there are $10$ lambencies for which the corresponding cases of umbral moonshine manifest a discriminant property, similar to that described for Thompson moonshine above. For $\ell=3$ for example, the corresponding umbral group is the unique non-trivial extension
\begin{gather}\label{eqn:G3}
1\to\ZZ/2\ZZ\to G^{(3)}\to M_{12}\to 1
\end{gather}
where $M_{12}$ is the sporadic simple {\em Mathieu group} with $\# M_{12}=95040$ (cf. \cite{atlas}). This double cover (\ref{eqn:G3}) of $M_{12}$ is typically denoted $2.M_{12}$. We have $m^{(3)}=3$, and the coefficients of $q^{\frac{11}{12}}$ in the McKay--Thompson series (\ref{eqn:umpmanalogue}) for $\ell=3$ realize the direct sum of two irreducible $16$-dimensional representations for $2.M_{12}$ whose character values lie in $\QQ(\sqrt{-11})$. Also, directly similar statements hold for $q^{\frac{8}{12}}$ and $q^{\frac{20}{12}}$, for the number fields $\QQ(\sqrt{-8})$ and $\QQ(\sqrt{-20})$, respectively.

So are there analogues of (\ref{eqn:jtau85995}) and (\ref{eqn:85995asatrace}), that could have served as starting points for umbral moonshine? Indeed 
there are. For example, we have
\begin{gather}\label{eqn:T316}
	T_3\left(\frac{1+\sqrt{-11}}6\right) = 9-16\left(\frac{-1+\sqrt{-11}}{2}\right)
\end{gather}
where $T_3$ is the normalized principal modulus for $\Gamma_0(3)$, as in (\ref{eqn:T3}), and the coefficient $C^{(3)}_{e,1}(-11)=2\times 16$ of $q^{\frac{11}{12}}$ in $H^{(3)}_{e,1}$ (\ref{eqn:vvHUM}) arises via twisted traces of singular moduli (cf. (\ref{eqn:trDf}), (\ref{eqn:trNDF})) as
\begin{gather}\label{eqn:16asatrace}
	16=\frac{1}{\sqrt{-11}}\left(T_3\left(\frac{-1+\sqrt{-11}}6\right) - T_3\left(\frac{1+\sqrt{-11}}6\right)\right).
\end{gather}
As the reader may guess, the CM points occurring as arguments of $T_3$ in (\ref{eqn:16asatrace}) are the $\tau_Q$ corresponding to representatives $Q$ of the two distinct orbits of $\Gamma_0(3)$ on positive-definite forms in $\mc{Q}_3(-11)$.

This relationship (\ref{eqn:16asatrace}) between umbral moonshine at $\ell=3$ and the principal modulus $T_3$ 
is just one consequence of the analysis in \cite{MR3357517}, wherein the complete realization of the mock modular forms $H^{(\ell)}_e$ 
via twisted traces of singular moduli of 
principal moduli for $\Gamma_0(\ell)$ 
is presented, for the so-called {\em pure A-type} cases, $\ell\in \{2,3,4,5,7,9,13,25\}$. The extension of this to arbitrary $\ell$ was subsequently explained in \cite{omjt}.

So from the point of singular moduli, Thompson moonshine and umbral moonshine are closely related. The key to getting from the former to the latter is to replace the $j$-function---a principal modulus for $\SL_2(\ZZ)$---with principal moduli for more general genus zero groups (cf. (\ref{eqn:XGamma})). 

Implicit here is the statement that every case of umbral moonshine has a genus zero group attached to it, and indeed, the notation in (\ref{eqn:lambencies}) was chosen 
with this in mind. For example, the lambencies of the form $\ell=m$ (being those in the first line of (\ref{eqn:lambencies})) correspond to the genus zero groups of the form $\Gamma_0(m)$. 
For a full description see \cite{mum}, wherein the association of genus zero groups to umbral moonshine was first described. 
For a conceptual explanation of 
this association see \cite{omjt}, wherein it is shown that genus zero groups naturally classify the optimal 
mock modular forms of weight $\frac12$ 
with rational Fourier coefficients. 

We refrain from explaining this classification in detail here, and refer to Theorems 1.2.1 and 1.2.2 of loc. cit. for the precise formulation. Instead we offer some remarks on the proof, since 
it 
involves notions we have encountered in \S\ref{sec:cuspformsellipticcurves}, including some that were pivotal for \S\ref{sec:ECONan}. 
For example, a key step 
is to show that the existence of an optimal mock modular form with shadow spanned by the theta functions $\theta^1_{m,r}$ of (\ref{eqn:theta1mr}) implies the vanishing of the $L$-function 
special values $L_f(1)$ (cf. (\ref{eqn:Lfs})) 
for all the cusp forms $f$ in a certain family. (This is Lemma 4.3.2 in loc. cit.) The proof depends on a certain precise formulation of Waldspurger's result (cf. the first paragraph of \S\ref{sec:cuspformsellipticcurves}) obtained by Skoruppa \cite{MR1074485}, following earlier work by Gross--Kohnen--Zagier \cite{MR909238} (cf. Theorem 3.3.3 in \cite{omjt}). With this in hand we are lead fairly directly (see Lemma 4.3.3 and 
Proposition 4.3.4 in \cite{omjt}) to a corresponding genus zero group. 

In the other direction, an optimal mock modular form $h$ whose shadow is not spanned by the $\theta^1_{m,r}$ has a cusp form $f$ of weight $2$ attached to it for which $L_f(1)\neq 0$. (This also depends on Skoruppa's formulation of Waldspurger's result.) By using this together with the main result of \cite{MR833192}---which also implies the infinitude of $E(\QQ)$ when $L_E(1)=0$ but $L'_E(1)\neq 0$ (cf. the last paragraph of \S\ref{sec:cuspformsellipticcurves})---we obtain, for some $m$, a $\Gamma_0(m)$-invariant divisor on $\HH$ which cannot occur as the divisor of a $\Gamma_0(m)$-invariant function on $\HH$. Then a result of Bruinier--Ono \cite{MR2726107} implies the existence of a transcendental Fourier coefficient for $h$. (See \S4.4 of \cite{omjt} for the details of this argument.)

We have seen now how the values of principal moduli at CM points may lead us to moonshine in half-integral weight, including Thompson moonshine, Mathieu moonshine and umbral moonshine in weight $\frac12$, and the more recently observed moonshine for O'Nan in weight $\frac32$. We have also seen how this connects to elliptic curve arithmetic in the latter case. So finally let us
consider the possibility that O'Nan moonshine reflects a 
more general relationship between finite groups and arithmetic invariants, as Mathieu moonshine manifests a special case of the more general umbral story. In contrast to the situation in weight $\frac12$, spaces of mock modular forms of weight $\frac32$ are expected to admit rational bases in general (see Theorem 4.7 of \cite{Bruinier201738} for a result in this direction), so there would be no genus zero classification of the kind we have described above. One may take the view that this makes O'Nan moonshine, and any weight $\frac32$ analogue of it, less special than the weight $\frac12$ counterparts. An alternative view is that a more general theory---wherein cusp forms make an impact via their presence rather than their absence---may prove useful to number theory and the study of elliptic curves.


\providecommand{\bysame}{\leavevmode\hbox to3em{\hrulefill}\thinspace}
\providecommand{\MR}{\relax\ifhmode\unskip\space\fi MR }
\providecommand{\MRhref}[2]{%
  \href{http://www.ams.org/mathscinet-getitem?mr=#1}{#2}
}
\providecommand{\href}[2]{#2}

\end{document}